\documentclass[12pt,psamsfonts]{amsart}
\usepackage{amsmath}
\usepackage{amsthm}
\usepackage{amssymb}
\usepackage{amscd}
\usepackage{amsfonts}
\usepackage{amsbsy}
\usepackage{epsfig,afterpage}
\usepackage{psfrag}

\begin{document}

\title[Animated phase portraits of nonlinear and chaotic dynamical systems]{Animated phase portraits of nonlinear and chaotic dynamical systems}

\author[J.M. Ginoux]
{Jean-Marc Ginoux$^1$}

\address{$^1$ Laboratoire {\sc Protee}, I.U.T. de Toulon,
Senior lecturer in Applied Mathematics,
Doctor in Mathematics, Doctor in History of Sciences,
Universit\'{e} du Sud, BP 20132, F-83957 La Garde cedex, France}
\email{ginoux@univ-tln.fr, http://ginoux.univ-tln.fr}.


\maketitle

\section{Introduction}

The aim of this section is to present programs allowing to highlight the \textit{slow-fast} evolution of the solutions of nonlinear and chaotic dynamical systems such as: Van der Pol, Chua and Lorenz models. These programs provide animated phase portraits in dimension two and three, i.e. ``integration step by step'' which are useful tools enabling to understand the dynamic of such systems.

\section{Van der Pol model}

The oscillator of B. Van der Pol \cite{VdP} is a second-order system with non-linear frictions which can be written:

\[
\frac{d^2y}{dt^2} - \alpha (1 - y^2)\frac{dy}{dt} + y = 0
\]

\smallskip

The particular form of the friction which can be carried out by an electric circuit causes a decrease of the amplitude of the great oscillations and an increase of the small. This equation constitutes the ``paradigm of relaxation-oscillations''. According to D'Alembert transformation \cite{alemb} any single $n^{th}$ order differential equation may be transformed into a system of $n$ simultaneous first-order equations and conversely. Let's consider that: $\alpha(1-y^2)\dot{y} = \frac{d}{dt} \alpha (y - \frac{y^3}{3})$ and let's pose: $x_1 = y$ and $y = -\alpha \dot{x_2}$. Thus, we have:

\[
\left[ \begin{array}{*{20}c}
 \dot{x_1} \vspace{4pt} \\
 \dot{x_2}
\end{array} \right] = \dfrac{1}{\alpha} \left[ \begin{array}{*{20}c}
 x_1 + x_2 - \dfrac{x_1^3}{3} \hfill \\
 - x_1
\end{array} \right]
\]

When $\alpha $ becomes very large, $x_1$ becomes a ``fast'' variable
and $x_2$ a ``slow'' variable. In order to analyze the limit $\alpha
\to \infty $, we introduce a small parameter $\varepsilon = 1
\mathord{\left/ {\vphantom {1 {\alpha ^2}}} \right.
\kern-\nulldelimiterspace} {\alpha ^2}$ and a ``slow time'' $t' = t
\mathord{\left/ {\vphantom {t {\alpha = }}} \right.
\kern-\nulldelimiterspace} {\alpha = }\sqrt \varepsilon t$. Thus,
the system can be written:

\begin{equation}
\label{eq1}
\left[ \begin{array}{*{20}c}
\varepsilon \dot{x_1}  \vspace{4pt} \hfill \\
\mbox{ } \dot{x_2} \hfill \\
\end{array} \right] = \left[ \begin{array}{*{20}c}
x_1 + x_2 - \dfrac{x_1^3}{3} \hfill \\
- x_1
\end{array} \right]
\end{equation}

\smallskip

with $\varepsilon$ a small positive real parameter $\varepsilon = 0.05$. System (\ref{eq1}) which has been extensively studied since nearly one century is called a \textit{slow-fast dynamical system} or a \textit{singularly perturbed dynamical system}\footnote{See for example Ginoux \cite{Gin}.}. Although, it has been established that system (\ref{eq1}) can not be integrated by quadratures (closed-form) it is well-known that it admits a solution of \textit{limit-cycle} type. The program presented here enables to emphasize the \textit{slow-fast} evolutions of the solution on this \textit{limit-cycle}.

\smallskip

\begin{center}
\begin{figure}[htbp]
\centerline{\includegraphics[width=9.8995cm,height=9.139cm]{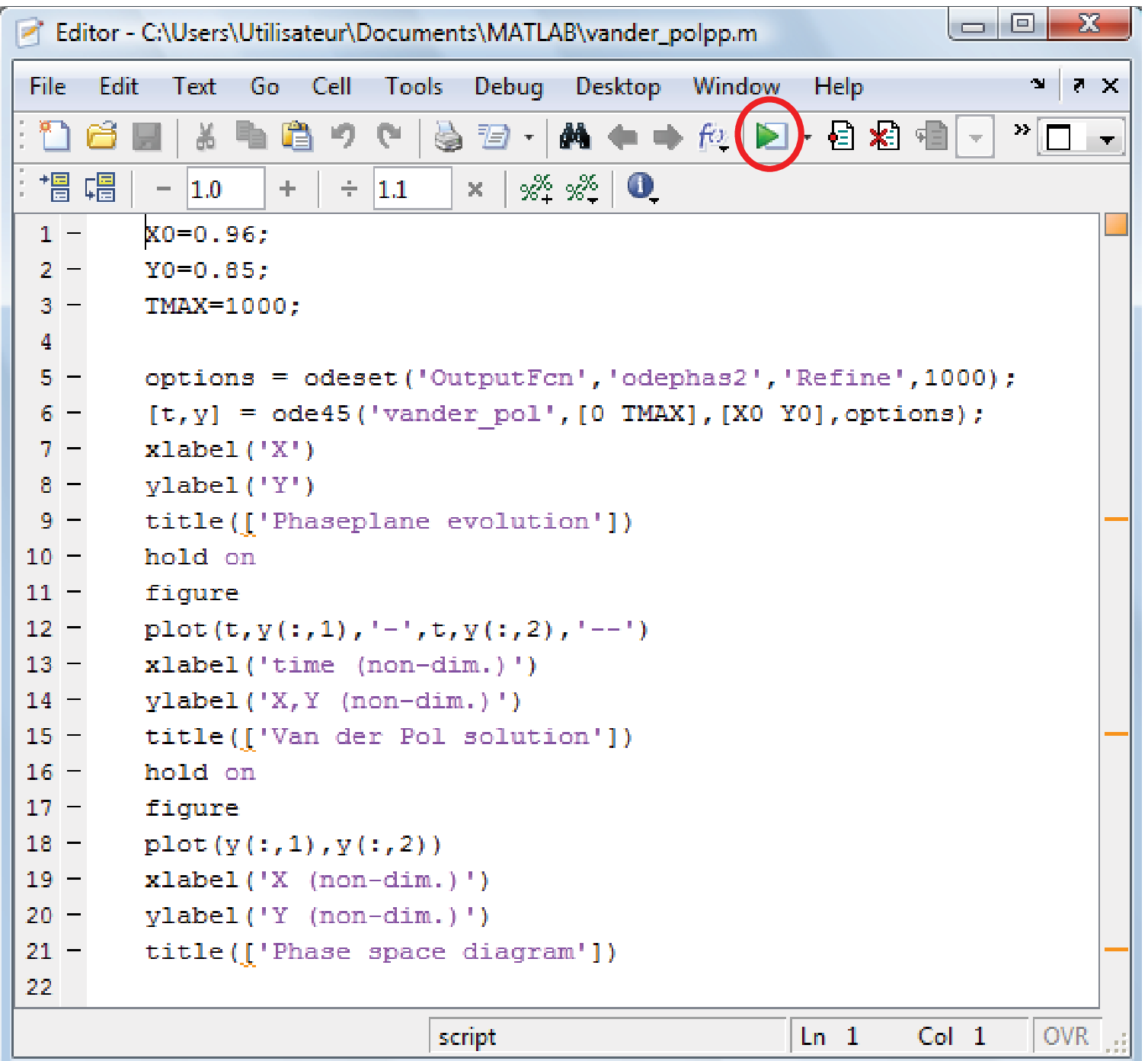}}
\label{fig1}
\caption{Program for 2D animated phase portrait}
\end{figure}
\end{center}

\begin{center}
\begin{figure}[htbp]
\centerline{\includegraphics[width=12.49cm,height=10.5cm]{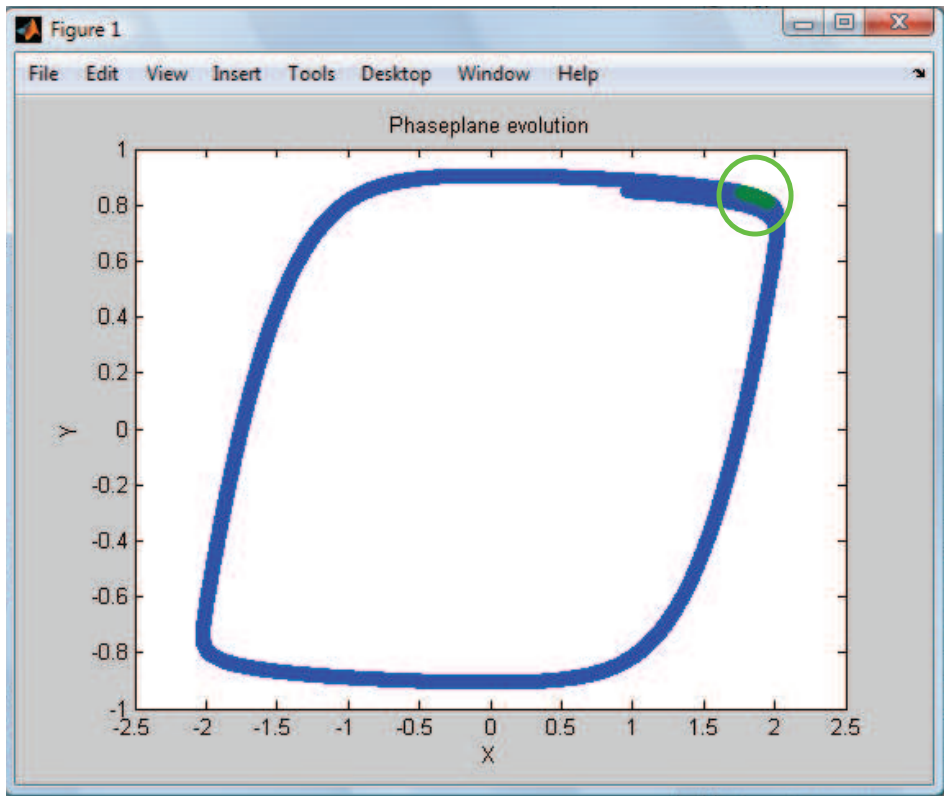}}
\label{fig2}
\caption{Animated phase portrait in 2D}
\end{figure}
\end{center}

First, copy the files named ``vanderpol'' and ``vanderpolpp'' into the ``current folder'' of MatLab. Then, open the m-file called ``vanderpolpp'' (see Fig. \ref{fig1}) and press the green button (red circle on the Fig. \ref{fig1}) to provide an animated plot 2D. On the Fig. \ref{fig2} the solution materialized by a green point (green circle on the Fig. \ref{fig2}) which evolves on the limit cycle, i.e., slowly on the nearly vertical parts and fast on the nearly horizontal parts.

All the main functions\footnote{See the textbook of Prof. Rudra Pratap \cite{Pratap}} are described in Tab. 1:

\begin{table}[htbp]
\begin{center}
\begin{tabular}{|p{5cm}|p{7cm}|}
\hline
\center{\textbf{Property}} &
\hspace{1.5cm} \textbf{Description}\\
\hline
\center{OutputFcn} &
A function for the solver to call after every successful integration step. \\
\hline
\center{odephas2} & 2D phase plane \\
\hline
\center{Refine} &
Increases the number of output points by a factor of Refine. \\
\hline
\end{tabular}
\label{tab1}
\caption{Description of the main functions}
\end{center}
\end{table}

The program provides the animated phase portrait corresponding to the solution of system (\ref{eq1}) and the time series.

\section{Chua's model}

The L.O. Chua's circuit \cite{Chua} is a relaxation oscillator with a cubic non-linear characteristic elaborated from a circuit comprising a harmonic oscillator of which operation is based on a field-effect transistor, coupled to a relaxation-oscillator composed of a tunnel diode. The modeling of the circuit uses a capacity which will
prevent from abrupt voltage drops and will make it possible to describe the fast motion of this oscillator by the following equations which also constitute a \textit{slow-fast dynamical system} or a \textit{singularly perturbed dynamical system}.

\begin{equation}
\label{eq2}
\left[ \begin{array}{*{20}c}
 \varepsilon \dot{x_1} \vspace{4pt} \hfill \\
 \mbox{ } \dot{x_2} \vspace{4pt} \hfill \\
 \mbox{ } \dot{x_3} \hfill \\
\end{array} \right] = \left[ \begin{array}{*{20}c}
  x_3 - \dfrac{44}{3}x_1^3 - \dfrac{41}{2}x_1^2 - \mu x_1 \vspace{4pt} \hfill \\
  - x_3 \vspace{4pt} \\
  - 0.7x_1 + x_2 + 0.24x_3 \\
\end{array} \right]
\end{equation}

\smallskip

with $\varepsilon$ and $\mu$ are real parameters $\varepsilon = 0.05$, $\mu = 2$. The system (\ref{eq2}) which can not be integrated by quadratures (closed-form) exhibits a solution evolving on ``chaotic attractor'' in the shape of a ``double-scroll''. The program presented here enables to emphasize the \textit{slow-fast} evolutions of the solution on the ``chaotic attractor''.

First, copy the files named ``chua'' and ``chua1'' into the ``current folder'' of MatLab. Then, open the m-file called ``chua1'' and press the green button to provide an animated plot 3D. The solution materialized by a green point evolves on the attractor according to slow and fast motion. The function ``odephas2'' is simply replaced by ``odephas3''.

\section{Lorenz model}

The purpose of the model established by Edward Lorenz \cite{Lorenz} was in the beginning to analyze the unpredictable behavior of weather. After having developed non-linear partial derivative equations starting from the thermal equation and Navier-Stokes equations, Lorenz truncated them to retain only three modes. The most widespread form of the Lorenz model is as follows:

\begin{equation}
\label{eq3}
\left[ \begin{array}{*{20}c}
 \dot{x_1} \hfill \\
 \dot{x_2} \hfill \\
 \dot{x_3} \hfill \\
\end{array} \right] = \left[ \begin{array}{*{20}c}
 \sigma \left( x_2 - x_1 \right) \\
  - x_1x_3 + rx_1 - x_2 \\
  x_1x_2 - \beta x_3 \\
\end{array} \right]
\end{equation}

\smallskip

with $\sigma$, $r$, and $\beta$ are real parameters: $\sigma = 10$, $\beta = \dfrac{8}{3}$, $r = 28$. 

\newpage

Although, this system is \underline{not} \textit{singularly perturbed} since it does not contain any small multiplicative parameter, it is a \textit{slow-fast dynamical system}. Its solution exhibits a solution evolving on ``chaotic attractor'' in the shape of a ``butterfly''. The program presented here enables to emphasize the \textit{slow-fast} evolutions of the solution on the ``chaotic attractor''.

\begin{table}[htbp]
\begin{center}
\begin{tabular}{|p{3cm}|p{9cm}|}
\hline
\center{\textbf{Name}} &
\hspace{1.5cm} \textbf{Description}\\
\hline
\center{vanderpol} &
Dynamical system (\ref{eq1}) with $\varepsilon = 1/20$. \\
\hline
\center{vanderpolpp} & animated phase portrait 2D \\
\hline
\center{chua} &
Dynamical system (\ref{eq2}) with $\varepsilon = 1/20$, $\mu = 2$. \\
\hline
\center{chua1} & animated phase portrait 3D \\
\hline
\center{lorenz} &
Dynamical system (\ref{eq3}) with $\sigma =10$, $\beta = 8/3$ $r = 28$ \\
\hline
\center{lorenz1} & animated phase portrait 3D \\
\hline
\end{tabular}
\label{tab2}
\caption{List of program for animated phase portraits.}
\end{center}
\end{table}

\end{document}